\documentclass{amsart}
\usepackage{amsmath,amscd,amssymb,latexsym,graphicx,color,xypic,mathrsfs}

\usepackage[english,francais]{babel}
\usepackage[latin1]{inputenc}

\def\intx{\overset{\:\circ}{X}}
\newcommand\tgt[1]{{}^{T}\kern-1pt #1}
\newcommand\adi[1]{{}^{ad}\kern-1pt #1}
\newtheorem{theorem}{Theorem}[section]

\newtheorem{e-proposition}[theorem]{Proposition}

\newtheorem{e-definition}[theorem]{Definition\rm}
\newtheorem{remark}{\it Remark\/}

\newtheorem{theoreme}{Th\'eor\`eme}[section]

\newtheorem{definition}[theoreme]{D\'efinition\rm}

\begin{document}

\def\to{\longrightarrow}

\def\L{\mathop{\wedge}}

\def\gpd{\,\lower1pt\hbox{$\longrightarrow$}\hskip-.24in\raise2pt
             \hbox{$\longrightarrow$}\,}

\begin{center}
{\Large\bf An index theorem for manifolds with boundary
}

\bigskip

\bigskip 

{\sc by Paulo Carrillo Rouse and Bertrand Monthubert}

\end{center}


\bigskip

\begin{center}
{\large \bf Abstract}
\end{center}
In \cite{Concg} II.5, Connes gives a proof of the Atiyah-Singer index theorem for closed manifolds by using deformation groupoids and appropiate actions of these on $\mathbb{R}^N$. Following these ideas, we prove an index theorem for manifolds with boundary. 


\vskip 0.5\baselineskip

\begin{center}
{\large \bf R\'esum\'e}
\end{center}
Dans \cite{Concg} II.5, Connes donne une preuve du th\'eor\`eme de l'indice d'Atiyah-Singer pour des vari\'et\'es ferm\'ees en utilisant des groupo\"ides de d\'eformation et des actions appropri\'ees de ceux-ci dans $\mathbb{R}^N$. Nous suivons ces id\'ees pour montrer un th\'eor\`eme d'indice pour des vari\'et\'es à bord.

\section*{Version fran\c{c}aise abr\'eg\'ee}
Dans \cite{Concg}, II.5, Alain Connes donna une preuve du théorème
d'Atiyah-Singer pour une variété fermée entièrement fondée sur l'utilisation de groupoïdes,
grâce à une action du groupoïde tangent de la variété sur $\mathbb{R}^N$. 
L'idée centrale est de remplacer des groupoïdes qui ne
sont pas (Morita) équivalents à des espaces, par des groupoïdes obtenus par
produit croisé et qui
possèdent cette propriété, ce qui permet ensuite de donner une
formule. 

Si $X$ est une variété à bord $\partial X$, nous construisons le
groupoïde $\mathcal{T}_bX:=(\adi{G_{\partial X}}\times
\mathbb{R})\bigcup_{\partial}TX$ en recollant   $\adi{G_{\partial
    X}}\times \mathbb{R}$ avec $TX$ le long de leur bord commun
$T\partial X\times \mathbb{R}$ (ici $\adi{G_{\partial
    X}}=T\partial X\cup \partial X\times \partial X\times (0,1)$ est le groupoïde
adiabatique).  Nous le recollons alors
avec le groupoïde tangent de l'intérieur de $X$,
$\tgt{G_{\intx}}=T\intx \cup \intx\times\intx\times (0,1]$ :
$\tgt{G_X}:=\mathcal{T}_bX\bigcup_{0}\tgt{G_{\intx}}$.

Il existe un homomorphisme
$\tgt{G_X}\stackrel{h}{\longrightarrow}\mathbb{R}^N$ induit par un
plongement de $X$ dans $ \mathbb{R}^{N-1}\times \mathbb{R}_+$, tel que
$\partial X$ se plonge dans $ \mathbb{R}^{N-1}\times
\mathbb{R}_+\times \{0\}$ et $\intx$ se plonge dans $ \mathbb{R}^{N-1}\times
\mathbb{R}^*_+$. Le produit croisé de $\tgt{G_X}$ par $h$ (noté $\tgt{(G_X)}_h$) est un
groupoïde propre dont les groupes d'isotropie sont triviaux, il est
donc Morita-équivalent à son espace d'orbites.

Soit $V(\intx)$ le fibré normal de  $\intx$ dans $\mathbb{R}^N$, et
$V(\partial X)$ le fibré normal de  $\partial X$ dans
$\mathbb{R}^{N-1}$ ; soit enfin  $V(X)=V(\intx)\bigcup V(\partial
X)$. En notant $\mathscr{D}_{\partial}=V(\partial X)\times
\{0\}\bigsqcup \mathbb{R}^{N-1}\times (0,1)$ et
$\mathscr{D}_{\circ}=V(\intx)\times \{0\}\bigsqcup \mathbb{R}^{N}\times
(0,1]$ les déformations au cône normal, on construit les espaces
$\mathscr{B}_{\partial}:=V(X)\bigcup_{\partial} \mathscr{D}_{\partial}$ et
$\mathscr{B}:=\mathscr{B}_{\partial}\bigcup_{\circ} \mathscr{D}_{\circ}$.

\begin{e-proposition}
  Le groupoïde $(\tgt{G_X})_h$ est Morita équivalent à l'espace $\mathscr{B}$.
\end{e-proposition}

Soit $$ind_f^X=(e_1)_*\circ(e_0)_{*}^{-1}:K^0(\mathcal{T}_bX)\longrightarrow K^0(\intx\times \intx)\approx \mathbb{Z}.$$
 \begin{definition}[Indice topologique pour une variété à bord]
Soit $X$ une variété à bord. L'indice topologique de $X$ est le morphisme
$$ind_t^X:K^0(\mathcal{T}_bX)\longrightarrow \mathbb{Z}$$
défini comme la  composition des trois morphismes suivants
\begin{enumerate}
\item L'isomorphisme de Connes-Thom  $CT_0$ suivi de l'équivalence de Morita  $\mathscr{M}_0$:
$$K^0(\mathcal{T}_bX)\stackrel{CT_0}{\longrightarrow}K^0((\mathcal{T}_bX)_{h_0})\stackrel{\mathscr{M}_0}{\longrightarrow}K^0(\mathscr{B}_{\partial}),$$
où $(\mathcal{T}_bX)_{h_0}$ est le produit croisé de $\mathcal{T}_bX$ par $h_0$ (l'homomorphisme h en $t=0$).
\item Le morphisme indice de l'espace de  déformation $\mathscr{B}$:
$
\xymatrix{
K^0(\mathscr{B}_{\partial})&
K^0(\mathscr{B})\ar[l]_-{(e_0)_*}^-{\approx}\ar[r]^-{(e_1)_*}& K^0(\mathbb{R}^N)
}
$
\item Le morphisme de périodicité de Bott :
$K^0(\mathbb{R}^N)\stackrel{Bott}{\longrightarrow}\mathbb{Z}.$
\end{enumerate}
\end{definition}

\begin{theorem}
  Pour toute variété à bord, on a l'égalité $$ind_f^X=ind_t^X.$$
\end{theorem}
\selectlanguage{english}

\section{Actions of $\mathbb{R}^N$}
\label{}
All the groupoids considered here will be  continuous family groupoids
\cite{LMN,paterson}. Hence we can consider their convolution and
$C^*$-algebras without any problem. If G is  such a groupoid, we will
denote by $K^0(G)$ the K-theory group of its $C^*$-algebra
(unless explicetely written otherwise). This is consistent with the usual notation when $G$ is a space
(a groupoid made only  of units). In the sequel, given a smooth
manifold $N$, we will denote by $\adi{G}_N:TN\times \{0\}\bigsqcup
N\times N\times \mathbb{R}^*\rightrightarrows N\times \mathbb{R}$, the
deformation to normal cone of $N$ in $N\times N$(for complete details about this 
deformation functor see \cite{Ca}). At each time, we will need to
restrict it to some interval, e.g. $[0,1]$ gives the tangent groupoid,
and $[0,1)$ gives the adiabatic groupoid.

Let $G\rightrightarrows M$ be a groupoid and 
$h:G\rightarrow \mathbb{R}^N$ a (smooth or continuous) homomorphism of
groupoids, ($\mathbb{R}^N$ as an additive group). Connes defined the
semi-direct product groupoid $G_h=G\times \mathbb{R}^N\rightrightarrows M\times \mathbb{R}^N$  (\cite{Concg}, II.5) with structure maps $t(\gamma,X)=(t(\gamma),X)$, $s(\gamma,X)=(s(\gamma),X+h(\gamma))$ and product 
$(\gamma,X)\circ (\eta,X+h(\gamma))=(\gamma\circ \eta,X)$ for composable arrows. 

At the level of $C^*$-algebras, $C^*(G_h)$ can be seen as the crossed product algebra 
$C^*(G)\rtimes \mathbb{R}^N$ where $\mathbb{R}^N$ acts on $C^*(G)$ by automorphisms by the formula: 
$\alpha_X(f)(\gamma)=e^{i\cdot (X\cdot h(\gamma))}f(\gamma)$, $\forall
f\in C_c(G)$, (see \cite{Concg}, propostion II.5.7 for details). In particular, in the case $N$ is even, we have a Connes-Thom isomorphism in K-theory
$K^0(G)\stackrel{\approx}{\rightarrow}K^0(G_h)$ (\cite{Concg}, II.C).

Using this groupoid, Connes gives a  conceptual, simple proof of the
Atiyah-Singer Index theorem for closed smooth manifolds. Let $M$ be a
smooth manifold, $G_M=M\times M$ its groupoid, and consider the tangent groupoid 
$\tgt{G_M}$. It is well known that the index morphism provided by this deformation groupoid is precisely the analytic index of Atiyah-Singer, \cite{Concg,MP}. In other words, the analytic index of $M$ is the morphism
\begin{equation}\label{index}
\xymatrix{
K^0(TM)\ar[r]^-{(e_0)_{*}^{-1}}&K^0(\tgt{G_M})\ar[r]^-{(e_1)_*}&K^0(M\times M)
=K^0(\mathscr{K}(L^2(M)))\approx \mathbb{Z},
}
\end{equation}
where $e_t$ are the obvious evaluation algebra morphisms at $t$. As discussed by
Connes, if the groupoids appearing in this interpretation of the index
were equivalent to spaces then we would immediately have a geometric
interpretation of the index. Now, $M\times M$ is equivalent to a point
(hence to a space), but the other fundamental groupoid playing a role
is not, that is, $TM$ is a groupoid whose fibers are the groups $T_xM$, which are  not equivalent (as groupoids) to a space. The idea of Connes is to use an appropriate action of the tangent groupoid in some $\mathbb{R}^N$ in order to translate the index (via a Thom isomorphism) in an index associated to a deformation groupoid which will be equivalent to some space. 

\section{Groupoids and Manifolds with boundary}
\label{bord}

Let $X$ be a manifold with boundary $\partial X$. We denote, as usual, $\intx$ the interior which is a smooth manifold. Let $X_{\partial}$ be the smooth manifold obtained by glueing $X$ with 
$\partial X\times [0,1)$ along their common boundary, $\partial X\sim \partial X\times \{0\}$. 
Set $TX:=TX_{\partial}|_{X}$, and consider the smooth manifold
$\mathcal{T}_bX:=(\adi{G_{\partial X}}\times
\mathbb{R})\bigcup_{\partial}TX$ obtained by glueing $\adi{G_{\partial
    X}}\times \mathbb{R}$ and $TX$ along their common boundary
$T\partial X\times \mathbb{R}$ ($\adi{G_{\partial X}}=T\partial X \cup
\partial X\times\partial X\times(0,1)$ is the adiabatic
groupoid). Now, we have a continuous family groupoid over $X_{\partial}$: 
$\mathcal{T}_bX
\rightrightarrows X_{\partial}$. As a groupoid it is the union of the  groupoids $\tgt{G_{\partial X}}\times \mathbb{R}\rightrightarrows \partial X\times [0,1)$ and $TX\rightrightarrows X$. For the topology, it is very easy to see that all the groupoid structures are compatible with the glueings we considered. 

We are going to consider a  deformation groupoid $\tgt{G_X}$ (\cite{Mont2}). This will be a natural generalisation of the Connes tangent groupoid of a smooth manifold, to the case with boundary. The space of arrows $\tgt{G_X}:=\mathcal{T}_bX\bigcup_{\circ}\tgt{G_{\intx}}$ is obtained by glueing at $T\intx$ ($T\intx\times\{0\}\subset \tgt{G_{\intx}}$ is closed). The space of units $X_{g_0}:=X_{\partial}\bigcup_{\circ}\intx\times [0,1]$ is obtained by glueing  $\intx\sim \intx\times \{0\}$ ($\intx\times \{0\}\subset \intx\times [0,1]$ is closed). Using the groupoid structures of $\mathcal{T}_bX\rightrightarrows X_{\partial}$ and of $\tgt{G_{\intx}}\rightrightarrows \intx\times [0,1]$, we have a continuous family groupoid $\tgt{G_X}\rightrightarrows X_{g_0}$. Again, all the groupoid structures are compatible with the considered glueings.

To define a homomorphism $\tgt{G_X}\stackrel{h}{\longrightarrow}\mathbb{R}^N$ we will need as in the nonboundary case an appropiate embedding. It is possible to find an embedding 
$i:X\hookrightarrow \mathbb{R}^{N-1}\times \mathbb{R}_+$ such that its restrictions to the interior 
and to the boundary are (smooth embeddings) of the following form 
$i_{\circ}:\intx\hookrightarrow \mathbb{R}^{N-1}\times \mathbb{R}_+^*$ and 
$i_{\partial}:\partial X\hookrightarrow \mathbb{R}^{N-1}\times \{0\}$.
We define the homomorphism $h:\tgt{G_X}\rightarrow \mathbb{R}^N$ as
follows.
\begin{equation}
  h:
  \begin{cases}
    h(x,X,0)=d_xi_{\circ}(X) \text{ and }
    h(x,y,\epsilon)=\frac{i_{\circ}(x)-i_{\circ}(y)}{\epsilon}  \text{ on }
    \tgt{G_{\intx}}\\
h(x,\xi,0,\lambda)=(d_xi_{\partial}(\xi),\lambda) \text{ and }
h(x,y,\epsilon,\lambda)=(\frac{i_{\partial}(x)-i_{\partial}(y)}{\epsilon},\lambda) \text{ on } \tgt{G_{\partial X}} \times
\mathbb{R} \\
h(x,X)=d_xi_{\circ}(X) \text{ on }  T\intx
  \end{cases}
\end{equation}

 Since all these morphisms are
compatible with the glueings, one has:

\begin{e-proposition}
With the formulas defined above, $h:\tgt{G_X}\rightarrow \mathbb{R}^N$ defines a homomorphism of continuous family groupoids.
\end{e-proposition}

The action of $\tgt{G_X}$ on $\mathbb{R}^N$ defined by $h$ is free because $i$ is an immersion. It is not necessarily proper (in the case of Connes \cite{Concg} II.5 it is since M was supposed closed), however we can prove the following: 

\begin{e-proposition}\label{propgrpd}
The groupoid $(\tgt{G_X})_h$ is a proper groupoid with trivial isotropy groups.
\end{e-proposition}

Notice that the groupoid $G_h$ is not the action groupoid (if not, the properness of the action would be equivalent to the properness of the groupoid). It is very important that the units of the groupoid $G_h$ be the units of $G$ times $\mathbb{R}^N$. 


As an immediate consequence of the propositions above, the groupoid
$(\tgt{G_X})_h$ is Morita equivalent to its space of orbits.  Let us specify this space. 

Let $V(\intx)$ be the total space of the normal bundle of $\intx$ in $\mathbb{R}^N$. 
Similarly, let $V(\partial X)$ be 
the total space of the normal bundle of $\partial X$ in $\mathbb{R}^{N-1}$. 
Observe that they have the same fiber vector dimension. In fact, their
union $V(X)=V(\intx)\bigcup V(\partial X)$, is a vector bundle over $X$, the normal bundle of $X$ 
in $\mathbb{R}^N$. 

Take $\mathscr{D}_{\partial}=V(\partial X)\times \{0\}\bigsqcup \mathbb{R}^{N-1}\times (0,1)$ the deformation to the normal cone associated to the embedding 
$\partial X\stackrel{i_{\partial}}{\hookrightarrow}  \mathbb{R}^{N-1}$.
We consider the space $\mathscr{B}_{\partial}:=V(X)\bigcup_{\partial}
\mathscr{D}_{\partial}$ glued over their common boundary  $V(\partial X)\sim V(\partial X)\times \{0\}$. On the other hand, take $\mathscr{D}_{\circ}=V(\intx)\times \{0\}\bigsqcup \mathbb{R}^{N}\times (0,1]$ the deformation to the normal cone associated to the embedding $\intx\stackrel{i_{\circ}}{\hookrightarrow}  \mathbb{R}^{N}$. We consider the space $\mathscr{B}:=\mathscr{B}_{\partial}\bigcup_{\circ} \mathscr{D}_{\circ}$ glued over  $V(\intx)$ by the identity map.

\begin{e-proposition}\label{espaceB}
The space of orbits of the groupoid $(\tgt{G_X})_h$ is $\mathscr{B}$. 
\end{e-proposition}

 We can give the explicit homeomorphism. The orbit space of
 $(\tgt{G_X})_h$ is a quotient of $X_{g_0}\times \mathbb{R}^N$. To define a map
 $\Psi:X_{g_0}\times \mathbb{R}^N \rightarrow \mathscr{B}$ it is
 enough to define it for each component of $X_{g_0}$. Let
 \begin{equation}
   \Psi:
   \begin{cases}
     \partial X \times (0,1) \times \mathbb{R}^{N-1}\times
 \mathbb{R}\rightarrow \mathbb{R}^{N-1}\times (0,1)\\
\Psi(a,t,\xi,\lambda):=(\frac{i_{\partial}(a)}{t}+\xi,t)
   \end{cases} 
   \begin{cases}
     \partial X \times \{0\} \times  \mathbb{R}^{N-1}\times
 \mathbb{R}\rightarrow
     V(\partial X)\\
\Psi(a,0,\xi,\lambda):=\overline{(i_{\partial}(a),\xi)}
   \end{cases} \end{equation}

 \begin{equation*}
    \begin{cases}
     \intx \times (0,1] \times \mathbb{R}^N\rightarrow
     \mathbb{R}^{N}\times (0,1]\\
\Psi(x,t,X):=(\frac{i_{\circ}(x)}{t}+X,t)
   \end{cases}
  \begin{cases}
     \intx \times \{0\} \times \mathbb{R}^{N}\rightarrow V(\intx)\\
\Psi(x,0,X):=\overline{(i_{\circ}(x),X)}
   \end{cases}
 \end{equation*}
where $\overline{\xi}$ denotes the class in 
$V_a(\partial
X):=\mathbb{R}^{N-1}/T_{i_{\partial}(a)}\partial
X$ (resp. $\overline{X}$ denotes the class in 
$V_x(\intx):=\mathbb{R}^{N}/T_{i_{\circ}(x)}\intx$). This gives a continuous map $\Psi:X_{g_0}\times \mathbb{R}^N \rightarrow \mathscr{B}$ that passes to the quotient into a homeomorphism $\overline{\Psi}:(X_{g_0}\times \mathbb{R}^N)/\sim \rightarrow \mathscr{B}$, where $(X_{g_0}\times \mathbb{R}^N)/\sim$ is the orbit space of the groupoid $(\tgt{G_X})_h$.

\section{The index theorem for manifolds with boundary}
\label{}

Deformation groupoids induce index morphisms. The groupoid $\tgt{G_X}$ is naturally parametrized by the closed interval $[0,1]$. Its algebra comes 
equipped with evaluations to the algebra of $\mathcal{T}_bM$ (at t=0) and to the algebra of $\intx\times \intx$ (for $t\neq 0$). We have a short exact sequence of $C^*$-algebras
\begin{equation}\label{btangentsuite}
\xymatrix{
0\ar[r]&C^*(\intx\times \intx\times (0,1])\ar[r]&C^*(\tgt{G_X})\ar[r]^-{e_0}&C^*(\mathcal{T}_bM)\ar[r]&0
}
\end{equation}
where the algebra $C^*(\intx\times \intx\times (0,1])$ is contractible. Hence applying the $K$-theory functor to this sequence we obtain an index morphism
$$ind_f^X=(e_1)_*\circ(e_0)_{*}^{-1}:K^0(\mathcal{T}_bX)\longrightarrow
K^0(\intx\times \intx)\approx \mathbb{Z}.$$

The morphism $h:\tgt{G_X}\rightarrow \mathbb{R}^N$ is by definition also parametrized by $[0,1]$, {\it i.e.}, we have morphisms $h_0:\mathcal{T}_bM\rightarrow \mathbb{R}^N$ and 
$h_t:\intx\times \intx \rightarrow \mathbb{R}^N$, for $t\neq 0$. We
can consider  the associated groupoids, which  satisfy the same
properties as in proposition \ref{propgrpd} (in fact, for proving such
proposition it is better to do it for each $t$, and to check all the compatibilities).

\begin{definition}\label{deftop}[Topological index morphism for a manifold with boundary]
Let $X$ be a manifold with boundary. The topological index morphism of $X$ is the morphism
$$ind_t^X:K^0(\mathcal{T}_bX)\longrightarrow \mathbb{Z}$$
defined (using an embedding as above) as the composition of the following three morphisms
\begin{enumerate}
\item The Connes-Thom isomorphism $CT_0$ followed by the Morita equivalence $\mathscr{M}_0$:
$$K^0(\mathcal{T}_bX)\stackrel{CT_0}{\longrightarrow}K^0((\mathcal{T}_bX)_{h_0})\stackrel{\mathscr{M}_0}{\longrightarrow}K^0(\mathscr{B}_{\partial})$$
\item The index morphism of the deformation space $\mathscr{B}$:
$
\xymatrix{
K^0(\mathscr{B}_{\partial})&
K^0(\mathscr{B})\ar[l]_-{(e_0)_*}^-{\approx}\ar[r]^-{(e_1)_*}& K^0(\mathbb{R}^N)
}
$\\
\item The usual Bott periodicity morphism:
$K^0(\mathbb{R}^N)\stackrel{Bott}{\longrightarrow}\mathbb{Z}.$\\
\end{enumerate}
\end{definition}

\begin{remark}
The topological index defined above is a natural generalisation of the topological index theorem defined by Atiyah-Singer. Indeed, in the boundaryless case, they coincide. 
The index of the deformation space $\mathscr{B}$ is quite easy to understand because we are dealing now with spaces (as groupoids the product is trivial), then the group $K^0(\mathscr{B})$ is the K-theory of the algebra of continuous functions vanishing at infinity $C_0(\mathscr{B})$ and the evaluation maps are completely explicit. In particular, if we identify $\mathscr{B}_{\partial}$ with an open subset of $\mathbb{R}^N$ (in the natural way), then the morphism $(ii)$ above correspond to the canonical extension of functions of $C_0(\mathscr{B}_{\partial})$ to $C_0(\mathbb{R}^N)$. 
\end{remark}

The following diagram, in which the morphisms $CT$ and $\mathscr{M}$ are the Connes-Thom and Morita isomorphisms respectively, is trivially commutative:

\begin{equation}\label{diagAPS}
\scriptsize{\xymatrix{
K^0(\mathcal{T}_bX) \ar[d]^-{\approx}_{CT}& K^0(\tgt{G_X}) \ar[d]^-{\approx}_{CT} \ar[l]_-{e_0}^-{\approx}\ar[r]^-{e_1}& 
K^0(\intx\times \intx)\ar[d]^-{\approx}_{CT}\\
K^0((\mathcal{T}_bX)_{h_0}) \ar[d]^-{\approx}_{\mathscr{M}}& K^0((\tgt{G_X})_h) \ar[d]^-{\approx}_{\mathscr{M}}\ar[l]_-{e_0}^-{\approx}\ar[r]^-{e_1}& K^0((\intx\times \intx))_{h_1})\ar[d]^-{\approx}_{\mathscr{M}}\\
K^0(\mathscr{B}_{\partial}) & K^0(\mathscr{B}) \ar[l]_-{e_0}^-{\approx}\ar[r]^-{e_1} & 
K^0(\mathbb{R}^N),
}}
\end{equation}

 The left vertical line gives the first part of the topological index
 map. The bottom line is the morphism induced by the
 deformation space $\mathscr{B}$. And the right vertical line is precisely the inverse of the Bott isomorphism $\mathbb{Z}=K^0(\{pt\})\approx K^0(\intx\times \intx)\rightarrow K^0(\mathbb{R}^N)$. Since the top
 line gives $ind_f^X$, we obtain the following result:
\begin{theorem}
For any manifold with boundary $X$, we have the equality of morphisms
$$ind_f^X=ind_t^X.$$
\end{theorem}

\section{Perspectives}
\label{}

As discussed in \cite{DL,DLN,LMN}, the index map $ind_f^X$ computes
the Fredholm index of a fully elliptic operator in the $b$-calculus of
Melrose.  We shall use the result proven here to give a formula in
relation to that of Atiyah-Patodi-Singer (\cite{APS}).

\end{document}